\documentclass[11pt,reqno]{amsart}

\usepackage{graphics}
\usepackage{amsmath,amsthm}
\usepackage{amscd}
\usepackage{amssymb}
\usepackage{euscript}

\numberwithin{equation}{subsection}

\newtheorem{thm}[subsection]{Theorem}

\newtheorem{prop}[subsection]{Proposition}
\newtheorem{cor}[subsection]{Corollary}

\theoremstyle{definition}
\newtheorem{definition}[subsection]{Definition}
\newtheorem{example}[subsection]{Example}
\newtheorem{remark}[subsection]{Remark}

\begin{document}

\title{Miniversal deformations of dialgebras}
\author{ Alice Fialowski and Anita Majumdar}
\begin{abstract}{We develop the theory of versal deformations of
dialgebras and describe a method for constructing a miniversal
deformation of a  dialgebra.}
\end{abstract}
\date{}

\address{Alice Fialowski, Institute of Mathematics, E$\ddot{o}$tv$\ddot{o}$s Lor$\acute{a}$nd University, 1117, Budapest, Hungary.
\\E-mail: {\tt fialowsk@cs.elte.hu}}
\address{Anita Majumdar, Dept. of Mathematics, Indian Institute of Science,
Bangalore-560012, India.\\E-mail: {\tt anita@math.iisc.ernet.in}}

\footnote{The first author was supported by grants OTKA T043641 and T043034.}
\footnote{The second author is supported by NBHM Post-doctoral fellowship.}

\maketitle{}
\section{Introduction}
The notion of Leibniz algebras and dialgebras was discovered by J.-L. Loday
while studying periodicity phenomena in algebraic K-theory \cite{L}. Leibniz
algebras are a non-commutative variation of Lie algebras and dialgebras are
a variation of associative algebras. Recall that any associative algebra
gives rise to a Lie algebra by $[x,y]=xy-yx$. The notion of dialgebras was
invented in order to build analogue of the couple
$$
\text{Lie algebras} \leftrightarrow \mbox{associative algebras},
$$
where Lie algebras are replaced by Leibniz algebras. Shortly, dialgebra is to Leibniz algebra, what associative algebra is to Lie algebra. A (co)homology
theory  associated to dialgebras was developed by J.-L. Loday, called
the dialgebra cohomology where planar binary trees play a crucial role
in the  construction. Dialgebra cohomology with coefficients was
studied by A.  Frabetti \cite{F1, F2} and deformations of dialgebras
were developed in \cite{MM1}. In the present paper, we develop a
deformation  theory of dialgebras over a commutative unital algebra
base, following \cite{Fi1},  and show that dialgebra cohomology is a
natural candidate for the  cohomology controlling the deformations. We
work out a construction of a versal deformation for dialgebras,
following \cite{FF}.

\medskip

The paper is organized as follows. In Section\,2, we recall some
facts on dialgebra and  its cohomology. In Section\,3, we
introduce the definitions of deformations  of dialgebras over a
commutative, unital algebra base. In Section\,4, we  produce an
example of an infinitesimal deformation of a dialgebra $D$ over a
field $K$ , denoted by $\eta_D$, and also show that this
deformation is co-universal in  the sense that, given any
infinitesimal deformation $\lambda$ of a dialgebra  $D$ with a
finite dimensional base $A$, there exists a unique homomorphism
$\phi: K\oplus HY^2(D,D)' \longrightarrow A$, where $HY^2(D,D)$
denotes the two dimensional cohomology of $D$ with coefficients
in itself, such that $\lambda$ is equivalent to the push-out
$\phi_*\eta_D.$ Section\,5 comprises results of Harrison
cohomology of a commutative unital algebra $A$ with coefficients
in a $A$-module $M$, which have been used in the paper. In
Section\,6 we introduce obstructions to extending a deformation
over a base $A$, to a deformation over a base $B$, where there
exists an extension $0\rightarrow K\stackrel{i}{\rightarrow} B
\stackrel{p}{\rightarrow} A\rightarrow 0$ of $A$. We show that an
obstruction  is a cohomology class, vanishing of which is a
necessary and sufficient  condition for the given deformation of
$D$ over base $A$ to be extended to a  deformation of $D$ over
base $B$. In Section\,7 we discuss extendible  deformations. Let
$\lambda$ be the deformation of $D$ over $A$ which is extendible.
We state that the two dimensional cohomology group $HY^2(D,D)$
operates transitively on the set of equivalence classes of
deformations $\mu$ of $D$  with base $B$ such that
$p_*\mu=\lambda.$ We also state that the group of  automorphisms
of the extension $0\rightarrow
K\stackrel{i}{\rightarrow}B\stackrel{p}{\rightarrow}A\rightarrow
0$ operates  on the set of equivalence classes of deformations
$\mu$ such that $p_*\mu=  \lambda$. These two actions are related
by a map called  the \textsl{differential} $d\lambda:
TA\rightarrow HY^2(D,D),$ where $TA$  denotes the tangent space
of $A$. In the last section, we present a construction of a
miniversal deformation of a dialgebra $D$.

\section{Dialgebra and its Cohomology}
Throughtout this paper, $K$ will denote the ground field of
characteristic zero. All tensor products shall be over $K$ unless
specified.  In this section, we recall the definition of a dialgebra
and the construction  of the dialgebra cochain complex. Since we are
interested in coefficients in  the dialgebra itself, we shall restrict
our definition to the same.

\begin{definition} A dialgebra $D$ over $K$ is a vector space over $K$
along with two  $K$-linear maps $\dashv : D\otimes D \longrightarrow D$
called left and  $\vdash : D\otimes D \longrightarrow D$ called right
satisfying the  following axioms :
\begin{equation}
\begin{array}{rcl}
x\dashv (y \dashv z)&\stackrel{1}=&(x\dashv y)\dashv z \stackrel{2}= x\dashv
(y \vdash z)\\
(x \vdash y) \dashv z &\stackrel{3}=& x \vdash (y \dashv z) \\
(x \dashv y)\vdash z &\stackrel{4}=&x \vdash (y \vdash z) \stackrel{5}=
(x\vdash y)\vdash z
\end{array}
\end{equation}
\\
for all $x,y,z\in D$.
\end{definition}

Apart from the known algebraic examples of dialgebras, \cite{L}, we
cite an interesting example of a family of dialgebras in the context of
functional analysis, \cite{RF}.
\begin{example}
Let ${\mathcal H}$ be a Hilbert space and $e \in \mathcal H$ with
$||e||=1$. Define two linear operations $\dashv$ and $\vdash$ by $$a
\dashv b = \left\langle b, e\right\rangle a, ~~~~~ a \vdash b =
\left\langle a, e\right\rangle b,$$ for $a, b \in {\mathcal H}$. Then
$({\mathcal H}, \dashv, \vdash)$ is a dialgebra, more precisely, a
normed dialgebra \cite{RF}.
\end{example}

A morphism $\phi : D \longrightarrow D'$ between two dialgebras is a
K-linear map such that $\phi(x\dashv y)=\phi(x) \dashv \phi(y)$ and
$\phi(x\vdash y)=\phi(x) \vdash \phi(y).$

A planar binary tree with $n$ vertices (in short, $n$-tree) is a planar
tree with $(n+1)$ leaves, one root and each vertex trivalent. Let $Y_n$
denote the set of all $n$-trees. Let $Y_0$ be the singleton set
consisting of a root only. The $n$-trees for $0\leq n\leq 3$ are given
by the following diagrams:

\hskip 1.9cm
\resizebox{19cm}{3.5cm}{\includegraphics[0in,0in][12in,2in]{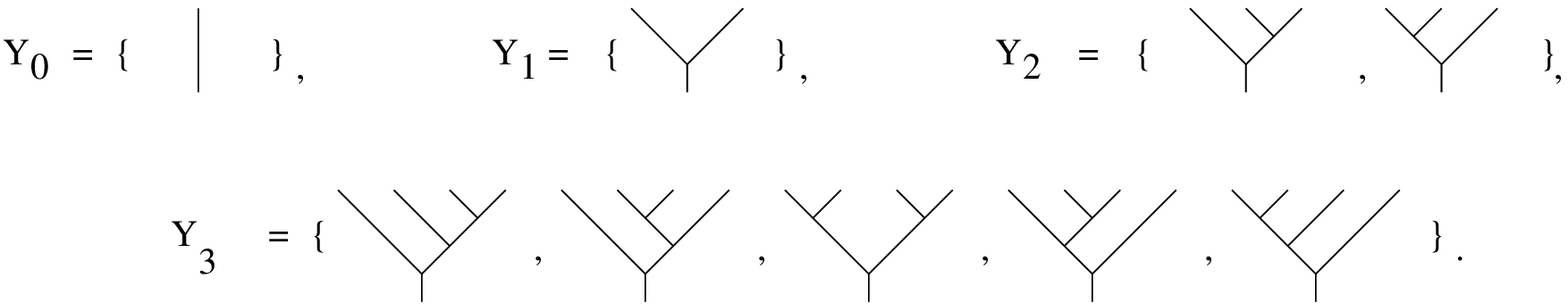}}

\bigskip
For any $y\in Y_n$, the $(n+1)$ leaves are labelled by $\{0,1,\ldots
,n\}$ from left to right and the vertices are labelled $\{1,2,\ldots
,n\}$ so that the $i$th vertex is between the leaves $(i-1)$ and $i$.
The only element $|$ of $Y_0$ is denoted by $[0]$ and the only element
of $Y_1$ is denoted by $[1]$. The grafting of a $p$-tree $y_1$ and a
$q$-tree $y_2$ is a $(p+q+1)$-tree denoted by $y_1\vee y_2$ which is
obtained  by joining the roots of $y_1$ and $y_2$ and creating a new
root from that vertex. This is denoted by $[y_1~ p+q+1~ y_2]$ with the
convention that all zeros are deleted except for the element in $Y_0$.
With this notation, the trees pictured above from left to right are
$[0],[1],[12],[21],[123],[213], [131], [312], [321]$.

For any $i$, $0\leq i\leq n$, there is a map, called the face map,  $d_i :
Y_n \longrightarrow Y_{n-1}$, $y\mapsto d_iy$ where $d_iy$ is obtained from
$y$ by deleting the $i$th leaf. The face maps satisfy the relations
$d_id_j=d_{j-1}d_i$, for all $i<j$.

\medskip

Let $D$ be a dialgebra over a field $K$. The cochain complex $CY^*(D,D)$
which defines the dialgebra cohomology $HY^*(D,D)$ is defined as follows.
For any $n\geq 0$, let $K[Y_n]$ denote the $K$-vector space spanned by $Y_n$
and $CY^n(D,D):= \mbox{Hom}_K(K[Y_n]\otimes D^{\otimes n},D)$ be the module
of $n$-cochains of $D$ with coefficients in $D$. The coboundary operator
$\delta : CY^n(D,D) \longrightarrow CY^{n+1}(D,D)$ is defined as the
$K$-linear map $\delta = \sum _{i=0}^{n+1}(-1)^i \delta^i$,
where \[(\delta^i f)(y;a_1,a_2,\ldots,a_{n+1})=\left\{\begin{array}{ll}
a_1 \circ_0^y f(d_0y;a_2,\ldots,a_{n+1}), & i=0\\
f(d_iy;a_1,\ldots,a_i \circ_i^y a_{i+1},\ldots,a_{n+1}),
& 1 \leq i\leq n\\
f(d_{n+1}y;a_1,\ldots,a_n)\circ_{n+1}^ya_{n+1}, & i=n+1
\end{array}
\right. \]
for any $y \in Y_{n+1}$; $a_1, \ldots,a_{n+1}\in D$
and $f:K[Y_n]\otimes D^{\otimes n} \longrightarrow D$.
Here, for any $i$, $0\leq i\leq n+1$. The maps $\circ_i:Y_{n+1}
\longrightarrow \{\dashv,\vdash\}$, are defined by

\[
\circ_0(y)= \circ_0^y:=\left \{\begin{array}{ll} \dashv & \mbox {if
$y$ is of the form $| \vee y_1$, for some $n$-tree $y_1$}\\ \vdash &
\mbox{otherwise} \end{array} \right.
\]

\[
\circ_i(y)= \circ_i^y:=\left \{\begin{array}{ll} \dashv & \mbox {if the
$i^{th}$ leaf of $y$ is oriented like `$\backslash $'}\\ \vdash &
\mbox{if the $i^{th}$ leaf of $y$ is oriented like `$/$' } \end{array}
\right.
\]
for $ 1\leq i\leq n$ and

\[ \circ_{n+1}(y)= \circ_{n+1}^y:=\left \{\begin{array}{ll}
\vdash & \mbox{if $y$ is of the form $y_1 \vee |$, for some $n$-tree
$y_1$}\\
\dashv & \mbox{otherwise}
\end{array},
\right. \]
where the symbol $`\vee'$ stands for grafting of trees \cite{L}.

There exists a pre-Lie algebra structure on $CY^*(D,D)$, \cite{MM1,
MM2}, the pre-Lie product being denoted by
$$
\circ: CY^n(D,D) \otimes CY^m(D,D)\longrightarrow CY^{n+m-1}.
$$
Also, if we modify the coboundary map $\delta$ by a sign, say
$dx=(-1)^{|x|}\delta(x)$, and define a bracket product on $CY^*(D,D)$
by $[x,y]= x\circ y -(-1)^{|x||y|} y\circ x$, which is the commutator
of the pre-Lie product, then $(CY^*(D,D), d)$ forms a differential
graded Lie algebra, \cite{MM2}, where $|x|= \mbox{deg}~x-1$.

\section{Deformations of Dialgebras}

Let $D$ be a dialgebra over $K$ and let $A$ be a commutative unital
algebra over $K$ with a fixed augmentation $\epsilon: A\longrightarrow
K$ with $\epsilon(1)=1.$ Let $\mathsf m = \ker\,\epsilon.$ We assume
$\dim\,({\mathsf m}^k/{\mathsf m}^{k+1})<\infty$, for all $k \geq 1$.

\begin{definition}
A deformation $\lambda$ of $D$ with base $(A,\mathsf m)$ is a dialgebra
structure on the tensor product $A\otimes_K D$ with the products
$\dashv_{\lambda}$ and $\vdash_{\lambda}$ being $A$-linear (or simply,
an $A$-dialgebra structure) such that $\epsilon\otimes \mbox{id}~:
A\otimes_K D\longrightarrow K\otimes D \cong D$ is a $A$-linear
dialgebra morphism. The left action of $A$ on $K\otimes D$ is given by
the augmentation map.
\end{definition}

We note that for $ x_1, x_2 \in D$, and $ a, b \in A$,
$$
a\otimes x_1 *_\lambda b\otimes x_2= ab ( 1\otimes x_1 * 1\otimes x_2),
$$
by A-linearity of the products, where $* =\{ \dashv, \vdash\}$. Also,
since $\epsilon \otimes ~\mbox{id}~: A\otimes D\longrightarrow K\otimes
D$ is a $A$-linear dialgebra homomorphism,
$$\begin{array}{rcl}
(\epsilon\otimes ~\mbox{id})~\{ 1\otimes x_1 *_\lambda 1\otimes x_2\}
&=&(\epsilon\otimes ~\mbox{id})(1\otimes x_1) * (\epsilon\otimes ~\mbox{id})(1\otimes x_2)\\
&=&(1\otimes x_1) * (1\otimes x_2)\\
&=&1\otimes (x_1* x_2)\\
&=&(\epsilon \otimes ~\mbox{id})(1\otimes (x_1* x_2)).
\end{array}
$$
So, $ (1\otimes x_1) *_\lambda (1\otimes x_2) - 1\otimes (x_1 * x_2)
\in \ker(\epsilon\otimes \mbox{id})$. Hence $ (1\otimes x_1) *_\lambda
(1\otimes x_2) = 1\otimes (x_1 * x_2)  + \sum_i m_i \otimes d_i$, where
$m_i \in \ker\,\epsilon= \mathsf m$ and $d_i \in D$ and $\sum_i m_i \otimes
d_i$ is a finite sum.

\begin{definition}
Two deformations of $D$ with the same base $A$ are called 
\textsl{equivalent} if there exists a $A$-linear dialgebra isomorphism
between the two copies of $A\otimes D$ with the two dialgebra
structures, compatible with $\epsilon \otimes \mbox{id}.$ A deformation
of $D$ with base $A$ is called \textsl{local} if the algebra $A$ is
local, and will be called \textsl{infinitesimal} if, in addition,
$\mathsf m^2=0$, where $\mathsf m$ is the maximal ideal of $A$.
\end{definition}

\begin{definition}
Let $A$ be a complete local algebra, that is, 
$A=
\overleftarrow{\displaystyle{\lim_{n\to\infty}}} (A/\mathsf m^n)$,
$\mathsf m$ denoting the maximal ideal in $A$. A formal deformation of $D$ with
base $A$ is a $A$-dialgebra structure on the completed tensor product
$A \widehat{\otimes}D=
\overleftarrow{\displaystyle{\lim_{n\to\infty}}} (A/\mathsf m^n)
\otimes D)$, such that $\epsilon
\widehat{\otimes}\mbox{id}~: A \widehat{\otimes} D\longrightarrow
K\otimes D=D$ is a $A$-linear dialgebra morphism.
\end{definition}

Two formal deformations of a dialgebra $D$ with the same base $A$ are called
\textsl{equivalent} if there exists a dialgebra isomorphism between the
two copies of $A \widehat{\otimes} D$ with the two dialgebra structures
compatible with $\epsilon \widehat{\otimes}\mbox{id}.$

\begin{example}
If $A= K[[t]]$ then a formal deformation of $D$ with base $A$ is the same as
a formal one-parameter deformation of $D$, \cite{MM1}.
\end{example}

Let $A'$ be a commutative algebra with identity, with a fixed augmentation
$\epsilon': A' \longrightarrow K$ and let $\phi: A\longrightarrow A'$ be an
algebra homomorphism, with $\phi(1)=1$ and $\epsilon' \circ \phi =
\epsilon.$ Then we can construct a deformation of $D$ with base $A'$ in the
following way.

\begin{definition} Let $\lambda$ be a deformation of the dialgebra $D$
with base $(A,\mathsf m)$. The \textsl{push}-\textsl{out} $\phi_*\lambda$
is the deformation of $D$ with base $(A',\mathsf m'=\ker\,\epsilon')$, which
is the dialgebra structure given by
\begin{align*}
a_1'\otimes_A(a_1\otimes x_1)
\dashv_{\phi_*\lambda}\ a_2'\otimes_A(a_2\otimes x_2)
=&\ a_1'a_2' \otimes_A (a_1\otimes x_1 \dashv_\lambda a_2\otimes x_2)\\
a_1'\otimes_A(a_1\otimes x_1)
\vdash_{\phi_*\lambda}\ a_2'\otimes_A(a_2\otimes x_2)
=&\ a_1'a_2' \otimes_A (a_1\otimes x_1 \vdash_\lambda a_2\otimes x_2),
\end{align*}
where $a_1, a_2\in A'$, $a_1, a_2 \in A$ and $l_1, l_2 \in D.$
Here we make use of the fact that $A'\otimes D= (A' \otimes_A A)\otimes D=
A'\otimes_A(A\otimes D),$ where $A'$ is regarded as an $A$-module by the
structure $a'a=a'\phi(a).$
\end{definition}

Similarly, one can define the \textsl{push}-\textsl{out} of
formal deformations.

\begin{remark}
We note that if the dialgebra structure $\lambda$ on $A\otimes D$
is given by
$$
(1\otimes x_1)*_\lambda (1\otimes x_2)= 1\otimes (x_1 * x_2)
 + \sum_{i=1}^n m_i \otimes d_i;~~ m_i \in\mathsf m, d_i \in D,
$$
then the dialgebra structure $\phi_*\lambda$ on $A'\otimes D$ is given by
$$
(1\otimes x_1)*_{\phi_*\lambda} (1\otimes x_2)
 = 1\otimes (x_1 * x_2)+ \sum_{i=1}^n \phi(m_i) \otimes d_i.
$$
\end{remark}

\bigskip

\section{Universal Infinitesimal and Miniversal Deformations of Dialgebras}

In \cite{FF}, the authors have produced a fundamental example of an
infinitesimal deformation of Lie algebras. Here we produce an example
of an infinitesimal deformation of dialgebras, which is obtained from
the aforesaid example, with slight modifications. Suppose $\dim\,HY^2(D,D) <\infty.$ This is, in particular, true if
$\dim\, D <\infty$. Consider the base of the deformation to be $A= K\oplus HY^2(D,D)',$ with $'$
denoting the linear dual. Here, $A$ is local with the maximal ideal $\mathsf m=
HY^2(D,D)',$ and $\mathsf m^2=0.$

Let
$$
\mu: HY^2(D,D) \longrightarrow CY^2(D,D)= \mbox{Hom\,}(K[Y_2]\otimes
D^{\otimes 2},D)
$$
which takes a cohomology class into a cocycle representing the class.
Define a dialgebra structure on
$$
\begin{array}{ll}
A\otimes D&= (K\oplus HY^2(D,D)')\otimes D\\
&=(K\otimes D) \oplus (HY^2(D,D)'\otimes D)\\
&=D \oplus (HY^2(D,D)'\otimes D)\\
&= D \oplus \mbox{Hom}~(HY^2(D,D),D)
\end{array}
$$
by

$$\begin{array}{ll}
(x_1,\phi_1)\dashv (x_2,\phi_2)&=(x_1\dashv x_2, \psi_{\ell})\\
(x_1,\phi_1)\vdash (x_2,\phi_2)&=(x_1\vdash x_2, \psi_r)
\end{array}
$$ where
$$
\begin{array}{ll}
\psi_{\ell}(\alpha)&=\mu(\alpha)([21]; x_1, x_2) + \phi_1(\alpha)\dashv x_2+
x_1\dashv \phi_2(\alpha)\\
\psi_r(\alpha)&=\mu(\alpha)([12]; x_1, x_2) + \phi_1(\alpha)\vdash x_2+
x_1\vdash \phi_2(\alpha),
\end{array}
$$
for $\alpha \in HY^2(D,D)$.

Using the dialgebra structure of $D$ and the fact that $\mu(\alpha)$ is a $2$-cocycle of $D$, one can check that the $\dashv$ and $\vdash$ products defined this way satisfy the
dialgebra axioms. It is to be noted that this deformation does not depend on the choice of
$\mu$, upto an isomorphism.

Let $\mu': HY^2(D,D) \longrightarrow CY^2(D,D)$ be another choice of $\mu.$
Define a homomorphism
$$\nu: HY^2(D,D) \longrightarrow CY^1(D,D)\cong \text{Hom}(D,D)$$ by
$\mu'(\alpha) - \mu (\alpha)= \delta\nu(\alpha),$ for all $\alpha \in
HY^2(D,D)$.
We define a linear automorphism $\rho$ of the space $A\otimes D= D\oplus
\text{Hom} (HY^2(D,D),D)$ by
$\rho(x, \phi)= (x, \psi)$ where $\psi(\alpha)= \phi(\alpha)+
\nu(\alpha)(x).$ It is straightforward to check that $\rho$ defines a
dialgebra isomorphism between the two dialgebra structures induced by $\mu$
and $\mu'$ respectively. We denote the infinitesimal deformation of $D$ as
constructed above by $\eta_D.$

\medskip

Below we will show the \textsl{couniversality} of $\eta_D$ in the
class of infinitesimal deformations:

Let $\lambda$ be an infinitesimal deformation of the dialgebra
$D$, with a finite dimensional local algebra base $A$, with
$\mathsf m^2=0$, where $\mathsf m$ is the maximal ideal of $A$.
Let $\xi \in {\mathsf m}'= \text{Hom}_K\,(\mathsf m,K).$ This is
equivalent to $\xi \in \text{Hom}_K~(A,K)$ with $\xi(1)=0.$

For $x_1, x_2 \in D,$ let us define a $2$-cochain as follows:
$$\alpha_{\lambda, \xi}([21]; x_1, x_2)= (\xi \otimes \mbox{id})((1\otimes
x_1)\dashv_{\lambda}(1\otimes x_2)) $$
and
$$\alpha_{\lambda, \xi}([12]; x_1, x_2)= (\xi \otimes \mbox{id})((1\otimes
x_1)\vdash_{\lambda}(1\otimes x_2)).$$

We claim that $\alpha_{\lambda,\xi}\in CY^2(D,D)$ is a $2$-cocycle. This is
because
\begin{align*}
\delta \alpha_{\lambda,\xi}([321];&x_1,x_2,x_3)\\
=\ &x_1\dashv
\alpha_{\lambda,\xi}([21];x_2,x_3)-\alpha_{\lambda,\xi}([21];x_1\dashv
x_2,x_3)\\
 &+ \alpha_{\lambda,\xi}([21]; x_1, x_2\dashv x_3)
-\alpha_{\lambda,\xi}([21]; x_1, x_2) \dashv x_3\\
=\ &x_1\dashv(\xi\otimes \mbox{id}~) ((1\otimes x_2)\vdash_{\lambda}
(1\otimes x_3))- (\xi \otimes \mbox{id})(1\otimes (x_1\dashv
x_2)\dashv_{\lambda}1\otimes x_3)\\
  &+ (\xi\otimes \mbox{id}) (1\otimes x_1 \dashv_{\lambda} 1\otimes
x_2\dashv x_3)- (\xi\otimes \mbox{id}~)(1\otimes x_1\dashv_{\lambda}1\otimes
x_2)\dashv x_3.
\end{align*}
If $\epsilon$ denotes the fixed augmentation of the algebra $A$, then
$$\epsilon \otimes \mbox{id}~:(1\otimes x_1 \dashv_{\lambda} 1\otimes x_2 -
1\otimes x_1\dashv x_2)=0,$$
i.e.\  $1\otimes x_1 \dashv_{\lambda} 1\otimes x_2 - 1\otimes x_1\dashv x_2
\in\mathsf m\otimes D.$
So,
\begin{align*}
(\xi\otimes \mbox{id}~) &((1\otimes x_1 \dashv_{\lambda} 1\otimes x_2)
\dashv_{\lambda} (1\otimes x_3))\\
=& (\xi\otimes \mbox{id}~) (((1\otimes x_1 \dashv x_2)+ \sum_i m_i\otimes
y_i)\dashv_{\lambda} (1\otimes x_3)))\\
=& (\xi\otimes \mbox{id}~) ((1\otimes x_1 \dashv
x_2)\dashv_{\lambda}(1\otimes x_3))+
(\xi\otimes \mbox{id}~) (\sum_i (m_i\otimes y_i)\dashv_{\lambda} (1\otimes
x_3))\\
=& (\xi\otimes \mbox{id}~) ((1\otimes x_1 \dashv
x_2)\dashv_{\lambda}(1\otimes x_3))+
(\xi\otimes \mbox{id}~) (\sum_i m_i(1\otimes y_i)\dashv_{\lambda} (1\otimes
x_3))\\
=&\alpha_{\lambda, \xi}([21];x_1\dashv x_2, x_3) + (\xi\otimes
\mbox{id})\sum_i m_i(1\otimes y_i \dashv_\lambda 1\otimes x_3).
\end{align*}
Note that in the second step from the end we make use of the action of
the algebra $A$ on $A\otimes D$.
\medskip

Now we have
$$1\otimes y_i \dashv_{\lambda} 1\otimes x_3 - 1\otimes y_i\dashv x_3 \in
m\otimes D,$$
$$1\otimes y_i \dashv_{\lambda} 1\otimes x_3 = 1\otimes y_i\dashv x_3 +h,$$
where $h \in \mathsf m\otimes D.$
Hence,
$$m_i(1\otimes y_i \dashv_{\lambda} 1\otimes x_3)= m_i (1\otimes y_i\dashv
x_3 +h).$$

Since $\mathsf m^2=0$, we have $m_i h=0.$ So, $m_i(1\otimes y_i \dashv_{\lambda}
1\otimes x_3)= m_i \otimes (y_i\dashv x_3),$ making use of the action of
$A$ on $A\otimes D.$
Next
\begin{align*}
(\xi\otimes \mbox{id}~)\sum_i m_i(1\otimes y_i\dashv_{\lambda} 1\otimes
x_3)&=\sum_i(\xi\otimes \mbox{id}~)(m_i\otimes y_i\dashv x_3)\\
&=\sum_i\xi(m_i)(y_i\dashv x_3)\\
&= \sum_i(\xi(m_i)y_i \dashv x_3)\\
&=(\xi\otimes \mbox{id}~)(\sum_im_i\otimes y_i)\dashv x_3\\
&=(\xi \otimes \mbox{id})\{(1\otimes x_1\dashv_\lambda 1\otimes
x_2)-1\otimes x_1\dashv x_2\}\dashv x_3\\
&=((\xi\otimes \mbox{id})(1\otimes x_1 \dashv_{\lambda}1\otimes x_2)\dashv
x_3\quad
\text{[using $\xi(1)=0$]}\\
&=\alpha_{\lambda, \xi}([12];x_1, x_2)\dashv x_3.
\end{align*}
Thus,
$$
\xi\otimes \mbox{id}~ ((1\otimes x_1 \dashv_{\lambda} 1\otimes x_2)
\dashv_{\lambda}(1\otimes x_3))
=\alpha_{\lambda, \xi}([21]; x_1\dashv x_2, x_3) + \alpha_{\lambda,
\xi}([21];x_1, x_2)\dashv x_3.
$$
In the same way,
$$
\xi\otimes \mbox{id}~ (1\otimes x_1 \dashv_{\lambda} (1\otimes
x_2\dashv_{\lambda}(1\otimes x_3))\\
=x_1 \dashv \alpha_{\lambda, \xi}([21]; x_2, x_3) + \alpha_{\lambda,
\xi}([21];x_1, x_2\dashv x_3).
$$
\medskip

 Since
 $$
 \xi\otimes \mbox{id}((1\otimes x_1 \dashv_{\lambda} 1\otimes
x_2)\dashv_{\lambda}(1\otimes x_3))-\xi\otimes\mbox{id}~(1\otimes x_1
\dashv_{\lambda} (1\otimes x_2\dashv_{\lambda}1\otimes x_3))=0,
$$
we have
$$
\delta\alpha_{\lambda, \xi}([321]; x_1, x_2, x_3)=0,
$$
and we can also show that $\delta\alpha_{\lambda, \xi}(y; x_1,
x_2, x_3)=0$ for all $y \in \{[312],[131], [213], [123]\}$.
\\
\\
The following proposition classifies all infinitesimal
deformations of $D$ over finite dimensional bases.
\begin{prop}
For any infinitesimal deformation $\lambda$ of a dialgebra $D$ with a finite
dimensional base $A$ there exists a unique homomorphism $\phi: K\oplus
HY^2(D,D)' \longrightarrow A$ such that $\lambda$ is equivalent to the
push-out $\phi_*\eta_D.$
\end{prop}

{\bf Proof.} Let $a_{\lambda, \xi} \in HY^2(D,D)$ be the cohomology class of
the cocycle
$\alpha_{\lambda, \xi}$, corresponding to $\xi\in\mathsf m'.$ Thus we have the
following homomorphisms:

$$\begin{array}{ll}
\alpha_\lambda:\mathsf m'&\longrightarrow CY^2(D,D)\\
a_\lambda:\mathsf m'&\longrightarrow HY^2(D,D).
\end{array}
$$

{\bf Step 1.} We show that the deformations $\lambda, \lambda'$
are equivalent if and only if $a_{\lambda}= a_{\lambda'}$. Let
$\lambda_1$ and $\lambda_2$ be two equivalent deformations of the
dialgebra $D$, with base $A$. By definition, there exists a
$A$-linear dialgebra isomorphism

\begin{equation}
\label{iso}
\rho: A\otimes D\longrightarrow A\otimes D,~ \text{such that}~
(\epsilon \otimes \mbox{id})\circ \rho= \epsilon \otimes \mbox{id}.
\end{equation}

Since $A\otimes D= D \oplus (\mathsf m\otimes D)$, the isomorphism  $\rho$ can
be written as $\rho = \rho_1 + \rho_2$ where $\rho_1: D\longrightarrow
D$ and $\rho_2: D \longrightarrow\mathsf m\otimes D$.

By using equation (\ref{iso}), we get $\rho_1=\mbox{id}.$ Note that by
the adjunction property of tensor products,
$$
\mbox{Hom}(D;\mathsf m\otimes D) \cong\mathsf m\otimes \mbox{Hom} (D,D)
\cong \mbox{Hom}(\mathsf m'; \mbox{Hom}(D,D)),
$$
where the isomorphisms are given by

\begin{equation}
\label{adiso}
\rho_2\longmapsto \sum_1^k m_i \otimes \phi_i\longmapsto \sum_i^k \chi_i.
\end{equation}
Here $\phi_i = (\xi_i \otimes \mbox{id})\circ \rho_2$ and $\chi_i
(\xi_j)= \delta_{i, j}\phi_i$, where $\{m_i\}_{1\leq i \leq k}$ is a
basis of $\mathsf m$ and $\{ \xi_j\}_{1\leq j\leq k}$ is a basis of
$\mathsf m'$.

We have by equation (\ref{adiso}),
$$\begin{array}{rcl}
\rho(1\otimes x)&=& \rho_1(1\otimes x) + \rho_2(1\otimes x)\\
&=& 1\otimes x + \sum_1^k m_i \otimes \phi_i(x).
\end{array}
$$

Using the notation $*=\{\dashv, \vdash\}$, the map $\rho$ is a
dialgebra homomorphism iff
$$
\rho( 1\otimes x_1 *_{\lambda_1} 1\otimes x_2)
= \rho( 1\otimes x_1) *_{\lambda_2} \rho( 1\otimes x_2),
$$

Let us set $\psi_i^r= \alpha_{\lambda_r, \xi_i}$, $i=1,2,\ldots, k$ and
$r=1,2.$ Then we have

\begin{equation}
1\otimes x_1 \dashv_{\lambda_r} 1\otimes x_2= 1\otimes x_1\dashv x_2 + \sum_i^k m_i \otimes \psi_i^r([21]; x_1, x_2)
\end{equation}
and
\begin{equation}
1\otimes x_1 \vdash_{\lambda_r} 1\otimes x_2= 1\otimes x_1\vdash x_2 + \sum_i^k m_i \otimes \psi_i^r([1 2]; x_1, x_2).
\end{equation}

Therefore, using the fact that $m_i. m_j=0$ for elements $m_i, m_j \in
\mathsf m$,
$$
\rho(1\otimes x_1 \dashv_{\lambda_1} 1\otimes x_2)
= 1\otimes x_1 \dashv x_2+ \sum_{i=1}^k m_i\otimes \phi_i(x_1 \dashv x_2)
+ \sum_{i=1}^k m_i ( 1\otimes \psi_i^1([21]; x_1, x_2)).
$$
\medskip
Similarly,
$$
\rho(1\otimes x_1 \vdash_{\lambda_2} 1\otimes x_2)
= 1\otimes x_1 \vdash x_2+ \sum_{i=1}^k m_i\otimes \phi_i(x_1 \vdash x_2)
+ \sum_{i=1}^k m_i ( 1\otimes \psi_i^1([12]; x_1, x_2)).
$$

Again,
\begin{align*}
\rho(1\otimes x_1)\dashv_{\lambda_2}&\ \rho(1\otimes x_2)\\
=& 1\otimes (x_1\dashv x_2) + \sum_{i=1}^k m_i
\otimes (\psi_i^2([21]; x_1, x_2)) + \sum_{i=1}^k m_i
\otimes \big(x_1\dashv \phi_i(x_2)\big)\\
& +\sum_{i=1}^k  m_i \otimes (\phi_i(x_1)\dashv x_2),
\end{align*}
and
\begin{align*}
\rho(1\otimes x_1)\vdash_{\lambda_2}\ & \rho(1\otimes x_2)\\
= & 1\otimes (x_1\vdash x_2) +
\sum_{i=1}^k m_i \otimes (\psi_i^2([21]; x_1, x_2)) +
 \sum_{i=1}^k m_i \otimes \big(x_1\vdash \phi_i(x_2)\big) \\
& +\sum_{i=1}^k  m_i \otimes (\phi_i(x_1)\vdash x_2).
\end{align*}

Thus, the following are equivalent:
\begin{align*}
\text{a)}&\quad \rho( 1\otimes x_1 \dashv_{\lambda_1} 1\otimes x_2)
 = \rho( 1\otimes x_1) \dashv_{\lambda_2} \rho( 1\otimes x_2)\\
\text{b)}&\quad \sum_{i=1}^k m_i \otimes (\psi_i^2([21]; x_1, x_2)
- \psi_i^1[21]; x_1, x_2)) +
 \sum_{i=1}^k m_i \otimes \delta \phi_i([21]; x_1, x_2)=0\\
\text{c)}&\quad \psi_i^1([21]; x-1, x_2) - \psi_i^2([21]; x_1, x_2) = \delta \phi_i ([21]; x_1, x_2)
\end{align*}

and similarly these are equivalent, too:
\begin{align*}
\text{a$'$)}&\quad \rho( 1\otimes x_1 \vdash_{\lambda_1} 1\otimes x_2)
= \rho( 1\otimes x_1) \vdash_{\lambda_2} \rho( 1\otimes x_2),\\
\text{b$'$)}&\quad \sum_{i=1}^k m_i \otimes (\psi_i^2([12]; x_1, x_2)
- \psi_i^1[12]; x_1, x_2))
 + \sum_{i=1}^k m_i \otimes \delta \phi_i([12]; x_1, x_2)=0\\
\text{c$'$)}&\quad \psi_i^1([12]; x-1, x_2) - \psi_i^2([12]; x_1, x_2) = \delta \phi_i ([12]; x_1, x_2).
\end{align*}

Hence,
$$\alpha_{\lambda_1, \xi_i}- \alpha_{\lambda_2, \xi_i}
= \delta \phi_i\quad \text{for $i\in\{1, 2,\ldots, k\}$}
\qquad \text{if and only if}\qquad
 a_{\lambda_1}= a_{\lambda_2}.
$$
This proves step 1.
\bigskip

{\bf Step 2}. Let $$\phi= \mbox{id} \oplus a_\lambda': K\oplus
HY^2(D,D)'\longrightarrow K\oplus\mathsf m= A.$$ Claim: $\phi_*
\eta_D$ is equivalent to $\lambda$. It follows from definitions
that $\alpha_{\phi_* \eta_D}= \mu\circ a_{\lambda}$. Thus,
$a_{\phi_* \eta_D}= a_\lambda$. Hence by step 1, $\phi_* \eta_D$
and $\lambda$ are isomorphic. This completes the proof of
Proposition 4.1.\qed

\medskip
Let $A$ be a local algebra with dim$(A/{\mathsf m}^2)<\infty.$ Then,
$A/{\mathsf m}^2$ is also local with the maximal ideal ${\mathsf
m}/{\mathsf m}^2$, and $({\mathsf m}/{\mathsf m}^2)^2=0.$
\begin{definition}
The linear dual space $\text{Hom}({\mathsf m}/{\mathsf m}^2, K)$ is
called the \textsl{tangent space} of $A$, and is denoted by $TA$.
\end{definition}

\begin{definition}\label{diff}
Let $\lambda$ be a deformation of $D$ with base $A$. Then the
mapping
$$a_{\pi_*\lambda}: TA=({\mathsf m}/{\mathsf m}^2)'\longrightarrow HY^2(D,D),$$
where $\pi$ is the projection $A\longrightarrow A/{\mathsf m}^2,$
is called the \textsl{differential} of $\lambda$ and is denoted
by $d\lambda.$
\end{definition}

\begin{definition}
A formal deformation $\eta$ of a dialgebra $D$ with base $B$ is called
\textsl{miniversal} if
\begin{enumerate}
\item for any formal deformation ${\lambda}$ of a dialgebra $D$ with any
local base $A$ there exists a homomorphism $f: B \rightarrow A$ such that
the deformation $\lambda$ is equivalent to $f_* \eta$;

\item with the above notations if $A$ satisfies the condition $\mathsf
m^2=0$, then $f$ is unique.
\end{enumerate}
If $\eta$ satisfies only condition (1), then it is called \textsl{versal}.
\end{definition}

The following proposition takes its shape from the general results
of Schlessinger \cite{Sch}. It was first shown for the case of Lie
algebras in \cite{Fi1}, and stated for Leibniz algebras in
\cite{FMM}. It is straightforward to see that it is true for the
case of dialgebras, too.

\begin{prop}
If the dimension of $HY^2(D,D)$ is finite, then there exists a
miniversal deformation of the dialgebra $D$.
\end{prop}\qed

\section{Some Facts about Harrison Cohomology}

Let $A$ denote a commutative algebra over $K$. In this section we
shall state a few results, without proof \cite{H}, about Harrison
cohomology groups of $A$ with coefficients in a $A$-module $M$.
Let $Ch(A)=\{Ch_q(A), \delta\}$ denote the Harrison complex of
$A$.
\begin{definition}
For an $A$-module $M$, the Harrison homology and cohomology of $A$
with coefficients in $M$ are defined as follows:
$$\begin{array}{ccc}
&H^{Harr}_q(A;M)&=H_q(Ch(A)\otimes M),\\
&H_{Harr}^q(A; M)&=H^q(\text{Hom}(Ch(A), M);
\end{array}
$$
\end{definition}

\smallskip

\begin{prop}
\begin{enumerate}
\item $H_{Harr}^1(A;M)$ is the space of derivations $A\rightarrow M$.

\item Elements of $H_{Harr}^2(A;M)$ correspond bijectively to isomorphism classes
of extensions $0\rightarrow M\rightarrow B\rightarrow A\rightarrow0$ of the
algebra $A$ by means of $M$.
\end{enumerate}
\end{prop}\qed

\begin{cor}\label{tan}
If $A$ is a local algebra with the maximal ideal $\mathsf m$, then
$$
H_{Harr}^1(A;K)=(\mathsf m/\mathsf m^2)'=TA.
$$
\end{cor}

\begin{prop}
Suppose $0\rightarrow M_r
\stackrel{\stackrel{-}{i}}{\rightarrow}B_{r-1}
\stackrel{\stackrel{-}{p}}{\rightarrow} A\rightarrow 0$ is an
$r$-dimensional extension of $A$. Then there is a
$(r-1)$-dimensional extension $0\rightarrow M_{r-1}
\stackrel{i}{\rightarrow}B_r \stackrel{p}{\rightarrow}
A\rightarrow 0$ of $A$ and a $1$-dimensional extension
$0\rightarrow K \stackrel{i'}{\rightarrow} B_r
\stackrel{p'}{\rightarrow}B_{r-1}\rightarrow 0$.
\end{prop}\qed

\begin{prop}
Let $0\rightarrow M\stackrel{i}{\rightarrow}B\stackrel{p}{\rightarrow}
A\rightarrow0$ be an extension of an algebra $A$ by $M$.
\begin{enumerate}
\item If $A$ has an identity
then so does $B$.
\item If $A$ is local with the maximal ideal $\mathsf m$, then $B$
is local with the maximal ideal $p^{-1}(\mathsf m).$
\end{enumerate}
\end{prop}\qed



\begin{definition}
Two extensions $B$ and $B'$ of the algebra $A$ by $M$ are said to be equivalent if there exists a $K$-algebra isomorphism $f: B\rightarrow B'$ such that the following diagram commutes.

$$\begin{array}{ccccccc}
0\longrightarrow & M & \stackrel{i_1}\longrightarrow & B &\stackrel{p_1}\longrightarrow & A & \longrightarrow 0\\
&&&&&&\\
 &\downarrow \mbox{id} & & \downarrow f & & \downarrow \mbox{id}&\\
&&&&&&\\
0\longrightarrow & M & \stackrel{i_2}\longrightarrow & B' &\stackrel{p_2}\longrightarrow & A & \longrightarrow 0.
\end{array}$$

An equivalence from $B$ to $B$ is said to be an automorphism of $B$ over $A$.
\end{definition}

\begin{prop}\label{aut}
$H_{Harr}^1(A;M)$ is isomorphic to the set of automorphisms of any given
extension $0\rightarrow M\stackrel{i}{\rightarrow}B\stackrel{p}{\rightarrow}
A\rightarrow0$ of $A$ by $M$.
\end{prop}\qed

\medskip


\section{Obstructions to Extending Deformations}

Let $A$ be a finite dimensional commutative, unital, local algebra with a fixed augmentation $\epsilon$, and let
$\lambda$ be a deformation of a dialgebra $D$ with base $A$. Let
$0\rightarrow K \stackrel{i}{\rightarrow}B \stackrel{p}{\rightarrow}A
\rightarrow 0$ be an extension of $A$, corresponding to a cohomology class
$f\in H^2_{Harr}(A;K)$. Let $q: A\rightarrow B$ be a splitting. Let $\widehat{\epsilon}: B\rightarrow K$ be the
augmentation of $B$. Let $I=i\otimes \mbox {id}~:
D=K\otimes D\rightarrow B\otimes D$ and $P=p\otimes \mbox {id}~:B\otimes
D\rightarrow A\otimes D$. Let $E= \widehat{\epsilon}\otimes \mbox {id}~:
B\otimes D\rightarrow K\otimes D=D$ and let $Q=q\otimes \mbox {id}~:
A\otimes D\rightarrow B\otimes D$. We define two $B$-bilinear operations
$\{~,~\}_{\dashv}$, $\{~,~\}_{\vdash}$ on $B\otimes D$ as follows:
\medskip

Let $l_1, l_2\in B\otimes D$.  Define
$$
\{l_1, l_2\}_\dashv= Q\{P(l_1)\dashv_\lambda P(l_2)\}
+I[I^{-1}(l_1-Q\circ P(l_1))\dashv I^{-1}(l_2- Q\circ P(l_2))],
$$
$$
\{l_1, l_2\}_\vdash= Q\{P(l_1)\vdash_\lambda P(l_2)\} +I[I^{-1}(l_1-Q\circ
P(l_1))\vdash I^{-1}(l_2- Q\circ P(l_2))].
$$
\smallskip
It is easy to verify that the two operations thus defined satisfy the
following properties:
\begin{align}\label{op}
&(i)\ P\{l_1, l_2\}_* = P(l_1)* P(l_2),
&\text{where\ }& *\in \{\dashv, \vdash\}, l_1, l_2\in B\otimes D,\\
&(ii)~~~\{I(l), l_1\}_* = I[l*E(l_1)],
&\text{where\ }& *\in \{\dashv, \vdash\}, l \in D, l_1\in B\otimes D.
\end{align}

Using the above two properties, one can show that
$$\begin{array}{rcl}
E\{l_1, l_2\}_{\dashv}&=& E(l_1) \dashv E(l_2)\\
E\{l_1, l_2\}_{\vdash}&=& E(l_1) \vdash E(l_2)\ .
\end{array}
$$

We define
\begin{align}\label{phi}
\phi([321];l_1, l_2, l_3)=& \{l_1, \{l_2, l_3\}_{\dashv}
\}_{\dashv}-\{\{l_1, l_2\}_{\dashv}, l_3\}_{\dashv},\\
\phi([312];l_1, l_2, l_3)=&\{\{l_1, l_2\}_{\dashv}, l_3\}_{\dashv}- \{l_1,
\{l_2, l_3\}_{\vdash} \}_{\dashv},\\
\phi([131];l_1, l_2, l_3)=& \{l_1, \{l_2, l_3\}_{\dashv}
\}_{\dashv}-\{\{l_1, l_2\}_{\vdash}, l_3\}_{\dashv},\\
\phi([213];l_1, l_2, l_3)=& \{\{l_1, l_2\}_{\dashv}, l_3\}_{\vdash}-\{l_1,
\{l_2, l_3\}_{\vdash} \}_{\vdash},\\
\phi([213];l_1, l_2, l_3)=&\{l_1, \{l_2, l_3\}_{\vdash}
\}_{\vdash}-\{\{l_1, l_2\}_{\vdash}, l_3\}_{\vdash}.
\end{align}

It is easy to see that $\phi(y; l_1, l_2, l_3) \in \ker\,P$ for all
$y  \in Y_3.$ Also, note that if any $l_i \in \ker\, E, i\in
\{1,2,3\}$, then  $\phi(l_1, l_2, l_3)=0.$ This defines the map
\begin{equation}\label{phibar}
\overline{\phi}: K[Y_3]\otimes D^{\otimes 3}=K[Y_3]\otimes ((B\otimes
D)/\ker\, E)^{\otimes 3}\rightarrow \ker\,P=D.
\end{equation}
Thus  $\overline{\phi}\in CY^3(D,D).$ One can check that $\delta
\overline{\phi}=0.$
\bigskip

Let $f'$ be cohomologous to $f$, and let $0\rightarrow K
\stackrel{i'}{\rightarrow}B' \stackrel{p'}{\rightarrow}A \rightarrow 0$ be
the extension corresponding to $f'$, which is isomorphic to the extension
corresponding to $f$. Since $B$ and $B'$ are isomorphic, without loss of
generality, we shall work with $B$.
\medskip

Let $\{~,~\}'_*, *\in \{\dashv, \vdash\}$ be another set of $B$-bilinear
operations on $B\otimes D$, satisfying ($1$) and ($2$) above. Then $\{l_1,
l_2\}'_*- \{l_1, l_2\}_* \in \ker\,P, *\in \{\dashv, \vdash\}$ for all
$l_1, l_2 \in B\otimes D.$ Also, $\{l_1, l_2\}'_*- \{l_1, l_2\}_* =0, *\in
\{\dashv, \vdash\}$ if $l_i \in \ker\,E, i\in \{1,2\}.$ This determines
a map $\psi: K[Y_2]\otimes D^{\otimes 2}=  K[Y_2]\otimes ((B\otimes
D)/\ker\, E)^{\otimes 2}\rightarrow \ker\,P=D.$ Thus, $\psi$
defines a $2$-cochain. Also, given an arbitrary $\psi\in CY^2(D,D)$, there
exists an appropriate $\{~,~\}_*'$ such that $\psi$ can be obtained as
$\{~,~\}_*'- \{~,~\}_*,$ where $*\in \{\dashv, \vdash\}.$

We remark here that if $\overline{\phi}$, $\overline{\phi}' \in CY^3(D,D)$
are the cochains corresponding to $\{~,~\}_*, \{~,~\}_*'$ in the sense of
the construction above, then
$$\overline{\phi}'-\overline{\phi}= \delta \psi.$$

Let ${\mathbb O}_\lambda (f) \in HY^3(D,D)$ be the cohomology class of the
cochain $\overline{\phi}.$ We define the following linear map.
$${\mathbb O}_\lambda : H^2_{Harr}(A,K)\longrightarrow HY^3(D,D),
~~~~~~~~~f\mapsto {\mathbb O}_\lambda(f).$$

We thus make the following proposition.

\begin{prop}
The deformation $\lambda$ with base $A$ can be extended to a deformation of
the dialgebra $D$ with base $B$ if and only if ${\mathbb O}_\lambda (f)=0.$
\end{prop} \qed

The cohomology class ${\mathbb O}_\lambda (f)$ is called the
\textsl{obstruction} to the extension of the deformation $\lambda$ from
$A$ to $B$.

\section{Extendible Deformations}
Let $A$ be a finite dimensional commutative, unital, local algebra with a fixed augmentation $\epsilon$, and let $\lambda$ be a deformation of a dialgebra $D$ with base $A$. Let $0\rightarrow K \stackrel{i}{\rightarrow}B \stackrel{p}{\rightarrow}A
\rightarrow 0$ be an extension of $A$, corresponding to a cohomology class
$f\in H^2_{Harr}(A;K)$.
Following the same arguments as in \cite{FF}, we can state the following
proposition.

\begin{prop}
$HY^2(D,D)$ operates transitively on the set of equivalence classes of
deformations $\mu$ of the dialgebra $D$ with base $B$ such that $p_*\mu=\lambda.$
\end{prop} \qed

We remark here that the group of automorphisms of the extension
$0\rightarrow
K\stackrel{i}{\rightarrow}B\stackrel{p}{\rightarrow}A\rightarrow 0$ is
$H_{Harr}^1(A;K),$ (\ref{aut}) and 
$H_{Harr}^1(A;K)=(\mathsf m/\mathsf m^2)'=TA,$ (\ref{tan}).
Note that by \ref{diff}, there exists a map $d\lambda: TA\rightarrow
HY^2(D,D).$ The group of automorphisms of the extension $0\rightarrow
K\stackrel{i}{\rightarrow}B\stackrel{p}{\rightarrow}A\rightarrow 0$ operates
on the set of equivalence classes of deformations $\mu$ such that $p_*\mu=
\lambda.$

We have the next proposition, the proof of which is straightforward.
\begin{prop}
The operation of $HY^2(D,D)$ on the set of equivalence classes of
deformations $\mu$ such that $p_*\mu=\lambda$ and the operation of the group
of automorphisms of the extension $0\rightarrow
K\stackrel{i}{\rightarrow}B\stackrel{p}{\rightarrow}A\rightarrow 0$ are
related by the differential $d\lambda: TA\rightarrow HY^2(D,D).$ In other
words, if $r: B\rightarrow B$ determines an automorphism of the extension
$0\rightarrow
K\stackrel{i}{\rightarrow}B\stackrel{p}{\rightarrow}A\rightarrow 0$ which
corresponds to an element $h\in H_{Harr}^1(A;K)=TA,$ then for any
deformation $\mu$ of $D$ with base $B$ such that $p_*\mu=\lambda,$ the
difference between the push-out $r_*\mu$ and $\mu$ is a cocycle of the
cohomology class $d\lambda(h).$
\end{prop} \qed

\begin{cor}
Suppose that the differential map $d\lambda: TA \longrightarrow HY^2(D,D)$ is onto. Then the group of automorphisms of
the extension $0\rightarrow
K\stackrel{i}{\rightarrow}B\stackrel{p}{\rightarrow}A\rightarrow 0$ operates
transitively on the set of equivalence classes of deformations $\mu$ of $D$
with base $B$ such that $p_*\mu=\lambda.$
\end{cor}

The proof of the following proposition is an imitation of the
proof presented in \cite{FMM}, for Leibniz algebras.

\begin{prop}
Let $A_1$ and $A_2$ be two finite dimensional local algebras with
augmentations $\epsilon_1$ and $\epsilon_2$, respectively. Let $\phi:
A_2 \longrightarrow A_1$ be an algebra homomorphism with $\phi(1)=1$
and $\epsilon_1 \circ \phi= \epsilon_2$. Suppose $\lambda_2$ is a
deformation of a dialgebra $D$ with base $A_2$ and $\lambda_1= \phi_*
\lambda_2$ is the push-out via $\phi$. Then the following diagram
commutes:
$$
\begin{array}{rcl}
H^2_{Harr}(A_1;K) & \stackrel{\phi^*}\longrightarrow & H^2_{Harr}(A_2;K)\\
& & \\
\theta_{\lambda_1}\searrow & & \swarrow\theta_{\lambda_2}\\
 & & \\
& HY^3(D;D) &
 \end{array} \ .
$$
\end{prop}
{\bf Proof.} Let $[f_{A_1}] \in H^2_{Harr}(A_1; K)$ correspond to the extension
$$0\longrightarrow K \stackrel{i_1}{\longrightarrow} A'_1 \stackrel{p_1}{\longrightarrow}A_1 \longrightarrow 0.
$$ 
Also, let $[f_{A_2}]= \phi^* ([f_{A_1}]) \in H^2_{Harr}(A_2; K)$ correspond to the extension $$0\longrightarrow K \stackrel{i_2}{\longrightarrow} A'_2 \stackrel{p_2}{\longrightarrow}A_2 \longrightarrow 0.$$ Let $q_k: A_k \longrightarrow A_k'$ be sections of $p_k$ for $k=1,2$. There exist $K$-module isomorphisms $A_k'\cong A_k \oplus K$. Let $(b,x)_{q_k}$ denote the inverse of $(b,x) \in A_k \oplus K$ under the isomorphisms. Define a linear map $\psi: A_2' \cong (A_2\oplus K) \longrightarrow A_1'\cong (A_1 \oplus K)$ by $\psi((a,x)_{q_2})= (\phi(a), x)_{q_1}$ for $(a,x)_{q_2}\in A_2'$. Thus we have a morphism of extensions

$$\begin{array}{ccccccc}
0\longrightarrow & K & \stackrel{i_2}\longrightarrow & A_2'&\stackrel{p_2}\longrightarrow & A_2 & \longrightarrow 0\\
&&&&&&\\
 &\downarrow \mbox{id}& & \downarrow \psi && \downarrow \phi &\\
&&&&&&\\
0\longrightarrow & K & \stackrel{i_1}\longrightarrow & A_1'&\stackrel{p_1}\longrightarrow & A_1 & \longrightarrow 0.
\end{array}$$

Let $I_k= i_k \otimes id$, $P_k= p_k\otimes id$ and $E_k=
\widehat{\epsilon}_k \otimes id$, where $\widehat{\epsilon}_k =
\epsilon_k \circ p_k$ for $k= 1,2$. If $m_{A_k}$ denote the unique
maximal ideal of $A_k$ then $m_{A_k'}= p_k^{-1}(m_{A_k})$ is the unique
maximal ideal of $A_k'$. Let the basis of $m_{A_k}$ and $m_{A_k'}$ be
$\{m_{k_i}\}_{1 \leq i\leq r_k}$ and $\{n_{k_i}\}_{1 \leq i\leq r_k+1}$
respectively, for $k=1,2$. Note that, $n_{k_j}= (m_{k_j}, 0)_{q_k}$ for
$1\leq j\leq r_k$ and $n_{k_{r_k+1}}= (0,1)_{q_k}$. The dialgebra
products on $A_2\otimes D$ is given by
$$\begin{array}{rcl}
(1\otimes x_1)\dashv_{\lambda_2} (1\otimes x_2)&=& 1\otimes (x_1 \dashv x_2) + \sum_{i=1}^{r_2} m_{2_i}\otimes \psi^2_i([21]; x_1, x_2)\\
(1\otimes x_1)\vdash_{\lambda_2} (1\otimes x_2)&=& 1\otimes (x_1 \vdash x_2) + \sum_{i=1}^{r_2} m_{2_i}\otimes \psi^2_i([12]; x_1, x_2)
\end{array}
$$ for $x_1, x_2 \in D$ and $\psi^2_i = \alpha_{{\lambda_2}, \xi_{2_i}}$, where $\{\xi_{2_i}\}$ is the dual basis of $\{m_{2_i}\}$.

Let $\phi(m_{2_i})= \sum_{j=1}^{r_1}c_{i,j}m_{1_j}$, $c_{i,j} \in K$
for $1\leq i\leq r_2$ and $1 \leq j \leq r_1$. Then the push-out
$\lambda_1= \phi_* \lambda_2$ on $A_1\otimes D$ is defined by
$$
\begin{array}{rcl}
(1\otimes x_1)\dashv_{\lambda_1} (1\otimes x_2)&=& 1\otimes (x_1 \dashv x_2) + \sum_{i=1}^{r_2}(\sum_{j=1}^{r_1} c_{i,j} m_{1_j})\otimes \psi^2_i([21]; x_1, x_2)\\
&=& 1\otimes (x_1 \dashv x_2) + \sum_{i=1}^{r_1} m_{1_j}\otimes \psi^1_j([21]; x_1, x_2)\\
(1\otimes x_1)\vdash_{\lambda_1} (1\otimes x_2)&=& 1\otimes (x_1 \vdash x_2) + \sum_{i=1}^{r_2}(\sum_{j=1}^{r_1} c_{i,j} m_{1_j})\otimes \psi^2_i([12]; x_1, x_2)\\
&=& 1\otimes (x_1 \vdash x_2) + \sum_{i=1}^{r_1} m_{1_j}\otimes \psi^1_j([12]; x_1, x_2)

\end{array}
$$
where $\psi^1_j\in CY^2(D,D)$ id defined by
$$\begin{array}{rcl}
\psi^1_j([21]; x_1, x_2)&=& \sum_{i=1}^{r_2} c_{i,j} \psi^2_i([21]; x_1, x_2)\\
\psi^1_j([12]; x_1, x_2)&=& \sum_{i=1}^{r_2} c_{i,j} \psi^2_i([12]; x_1, x_2)
\end{array}
$$ for $x_1, x_2 \in D$.

For any $2$-cochain $\chi \in CY^2(D,D)$, let us define $A_k'$ bilinear operations $\{, \}_{\dashv, k}, \{, \}_{\vdash, k}: (A_k'\otimes D)^{\otimes 2}\longrightarrow A_k' \otimes D$ by lifting $\lambda_k$,
$$
\begin{array}{rcl}
(1\otimes x_1)\dashv_{k} (1\otimes x_2)&=& 1\otimes (x_1 \dashv x_2) + \sum_{j=1}^{r_k} n_{k_j}\otimes \psi^k_j([21]; x_1, x_2)\\
&+& n_{k_{r_k+1}} \chi([21]; x_1, x_2)\\
(1\otimes x_1)\vdash_{k} (1\otimes x_2)&=& 1\otimes (x_1 \vdash x_2) + \sum_{j=1}^{r_k} n_{k_j}\otimes \psi^k_j([12]; x_1, x_2)\\
&+& n_{k_{r_k+1}} \chi([12]; x_1, x_2)
\end{array}
$$
for $k=1,2$ and $x_1, x_2 \in D$. The operations $\{, \}_{\dashv, k}, \{, \}_{\vdash, k}$, for $k=1,2$ satisfy the conditions (i) and (ii) of \ref{op}.

We shall show that
$\psi \otimes id: A_2' \otimes D \longrightarrow A_1' \otimes D$
preserves the liftings. It is enough to show that
$$(\psi \otimes id)(1\otimes x_1 *_2 1\otimes x_2)
= \psi \otimes id(1\otimes x_1) *_1 \psi \otimes id(1\otimes x_2),$$
where $* \in \{\dashv, \vdash\}$ and $x_1, x_2 \in D$.

Now
\begin{align*}
(\psi \otimes id)(&1\otimes x_1 \dashv_2 1\otimes x_2)
 =\psi(1)\otimes (x_1 \dashv x_2) \\
 &+ \sum_{j=1}^{r_2}\psi(1) \psi(n_{2_j})\otimes \psi_j^2([21]; x_1, x_2)\
+ \psi(1) \psi(n_{2_{r_2+1}}) \otimes \chi([21]; x_1, x_2)\\
=&1\otimes (x_1 \dashv x_2)+ \sum_{j=1}^{r_2}
\Big(\sum_{i=1}^{r_1} c_{j,i} m_{1_i}\Big) \otimes \psi_j^2([21]; x_1,
x_2) + n_{1_{r_1+1}}\otimes \chi([21]; x_1, x_1),
\end{align*}
where we used that $\phi(m_{2_j})= \sum_{i=1}^{r_1} c_{j,i} m_{1_i}$
and
$$
\psi(n_{2_{r_2+1}})= \psi((0,1)_{q_2})= (\phi(0), 1)_{q_1}
= n_{1_{r_1+1}}.
$$
Simplifying, we conclude that
\begin{align*}
(\psi \otimes& id)(1\otimes x_1 \dashv_2 1\otimes x_2)\\
 &=\psi(1) \otimes (x_1 \dashv x_2)+
 \sum_{i=1}^{r_1} \psi(1) m_{1_i} \otimes \psi_i^1([21]; x_1, x_2)
+ \psi(1) n_{1_{r_1+1}}\otimes \chi([21]; x_1, x_2)\\
&=(\psi(1)\otimes x_1) \dashv_1 (\psi(1)\otimes x_2) \\
&=\psi \otimes id(1\otimes x_1) \dashv_1 \psi \otimes id(1\otimes x_2).
\end{align*}

Similarly, we can show that
$$
(\psi \otimes id)(1\otimes x_1 \vdash_2 1\otimes x_2)= \psi
\otimes id(1\otimes x_1) \vdash_1 \psi \otimes id(1\otimes x_2).
$$
Let $\phi_k$ be defined by the operations $\{, \}_{\dashv, k}, \{,
\}_{\vdash, k}$ as has been defined in \ref{phi} and
$\overline{\phi_k}$ the corresponding cocycle as in \ref{phibar}.
Since, $\psi(n_{2_{(r_2+1)}})= n_{1_{(r_1+1)}}$, we have
$[\overline{\phi_2}]= [\overline{\phi_1}]$. Therefore,
$$
\theta_{\lambda_1}([f_{A_1}])=[\overline{\phi_1}] =
[\overline{\phi_2}]= \theta_{\lambda_2}([f_{A_2}])=
\theta_{\lambda_2}\circ \phi^*([f_{A_1}]).
$$
Hence, $\theta_{\lambda_1}=\theta_{\lambda_2}\circ \phi^*$. \qed

\section{Construction of a Miniversal Deformation of a Dialgebra}
An explicit description of the construction of a versal deformation of
a Lie algebra is given in \cite{FF}, and of a Leibniz algebra is given
in \cite{FMM}. Here we sketch the construction, for the case of a
dialgebra, following the same techniques developed in \cite{FF},
\cite{FMM}.

\medskip
Start with a dialgebra $D$ with $\mbox{dim}(HY^2(D,D)) < \infty$.
Consider the extension
$$
0\longrightarrow HY^2(D,D)' \stackrel{i}{\longrightarrow} C_1 \stackrel{p}
{\longrightarrow} C_0 \longrightarrow 0,
$$
where $C_0=K$, $C_1= K\oplus HY^2(D,D)'$. Let $\eta_1$ denote the
universal infinitesimal deformation with base $C_1$ as described in
Section~4.

\medskip
Suppose for some $k\geq 1$, we have constructed a finite
dimensional local algebra $C_k$, and a deformation $\eta_k$ of
$D$ with base $C_k$. Let
$$
\mu : H^2_{\text{\it Harr}}(C_k; K) \longrightarrow (Ch_2(C_k))'
$$
be a homomorphism mapping a cohomology class into a cocycle
representing the class. The dual map of $\mu$
$$
f_{C_k}: Ch_2(C_k)\longrightarrow  H^2_{\text{\it Harr}}(C_k; K)'
$$
corresponds to the following extension of $C_k$:
$$
0\longrightarrow H^2_{\text{\it Harr}}(C_k; K)'
\stackrel{\overline{i}_{k+1}}{\longrightarrow} \overline{C}_{k+1}
\stackrel{\overline{p}_{k+1}}{\longrightarrow} C_k \longrightarrow 0.
$$
The obstruction $\theta([f_{C_k}])\in H^2_{\text{\it Harr}}(C_k;
K)'\otimes HY^3(D,D)$ yields a map $\omega_k: H^2_{\text{\it Harr}}(C_k;
K)\longrightarrow HY^3(D,D)$,  by adjunction property of tensor
products, with the dual map
$$
\omega_k': HY^3(D,D)' \longrightarrow H^2_{\text{\it Harr}}(C_k; K)'.
$$
This induces the following extension
$$
0\longrightarrow \text{coker}(\omega_k')\longrightarrow \overline{C}_{k+1}/\overline{i}_{k+1}\circ \omega_k'(HY^3(D,D)')\longrightarrow C_k \longrightarrow 0.
$$
This yields the extension
$$
0\longrightarrow (\ker(\omega_k))' \stackrel{i_{k+1}}{\longrightarrow}C_{k+1}\stackrel{p_{k+1}}{\longrightarrow}C_k \longrightarrow 0$$ where $C_{k+1}= \overline{C}_{k+1}/ \overline{i}_{k+1}\circ \omega_k'(HY^3(D,D)')$ and $i_{k+1}, p_{k+1}$ are the mappings induced by $\overline{i}_{k+1}$ and $\overline{p}_{k+1}$ respectively. Along the same lines as in \cite{FF}, \cite{FMM} we have the following proposition:
\begin{prop}
The deformation $\eta_k$ with base $C_k$ of a dialgebra $D$ admits an
extension to a deformation with base $C_{k+1}$ which is unique up to an
isomorphism and an automorphism of the extension
$$
0\longrightarrow (\ker(\omega_k))'\stackrel{i_{k+1}}{\longrightarrow}C_{k+1}\stackrel{p_{k+1}}{\longrightarrow}C_k \longrightarrow 0.$$
\end{prop}\qed

This process gives rise to a sequence of finite dimensional local
algebras $C_k$ and deformations $\eta_k$ of the dialgebra $D$ with base
$C_k$
$$
K \stackrel{p_1}{\longleftarrow} C_1 \stackrel{p_2}{\longleftarrow} C_2 \stackrel{p_3}{\longleftarrow} \ldots \stackrel{p_k}{\longleftarrow} C_k \stackrel{p_{k+1}}{\longleftarrow} C_{k+1} \ldots
$$
such that $p_{k+1}*\eta_{k+1}= \eta_k$. By taking the projective limit
we obtain a formal deformation $\eta$ of $D$ with base
$C=\overleftarrow{\displaystyle{\lim_{k\to\infty}}} C_k$.

\medskip

Let $\dim\,(HY^2(D,D))=n$ and $K[[HY^2(D,D)']]$ denote the formal
power series ring in $n$ variables.  Also let $m$ denote the unique
maximal ideal in $K[[HY^2(D,D)']]$, consisting of all elements with
constant term zero. We have the following proposition, whose proof  can
be found in \cite{FMM}.
\begin{prop}
The complete local algebra
$C= \overleftarrow{\displaystyle{\lim_{k\to\infty}}} C_k$
can be described as
$$
C\cong K[[HY^2(D,D)']]/I,
$$
where $I$ is an ideal contained in $m^2$.
\end{prop}
\qed

Along the same lines as in \cite{FF}, \cite{FMM}, we state the
following theorem, proof of which obeys the same techniques as
developed in \cite{FF}.
\begin{thm}
Let $D$ be a dialgebra with $\dim(HY^2(D,D))< \infty$. Then the
formal deformation $\eta$ with base $C$ as described above is a
miniversal deformation of $D$.
\end{thm}\qed

\medskip

\noindent{\sc Acknowledgements}

\smallskip
We would like to thank Professor Goutam Mukherjee for his
interest in our work and the Indian Statistical Institute,
Calcutta, for hospitality.

\bigbreak

\end{document}